\documentclass[reqno]{amsart}
\usepackage{amsthm}
\usepackage{epsfig}
\usepackage{pst-grad} 
\usepackage{pst-plot} 
\usepackage[english]{babel}
\usepackage[latin1]{inputenc}
\usepackage{indentfirst}
\usepackage{graphicx,subfigure}
\usepackage{tabularx}
\usepackage{xspace} 
\usepackage{amsmath,amssymb}
\usepackage{color}
\definecolor{oneblue}{rgb}{0,0.0,0.75}

\usepackage{fancyhdr}
 \usepackage{setspace}

\makeatletter
\newcommand{\sech}{\mathop{\operator@font sech}}
\newcommand{\sign}{\mathop{\operator@font sign}}
\makeatother

\newtheorem{lemma}{Lemma}[section]
\newtheorem{theorem}{Theorem}[section]

\newtheorem{remark}{Remark}[section]
\numberwithin{equation}{section}

\begin{document}

\title[Semidiscrete approximation to Benjamin-type equations]{Semidiscrete approximation to Benjamin-type equations}

\author[V. A. Dougalis]{Vassilios A. Dougalis}
\address{Mathematics Department, University of Athens, 15784
Zographou, Greece \and Institute of Applied \& Computational
Mathematics, FO.R.T.H., 71110 Heraklion, Greece}
\email{doug@math.uoa.gr}

\author[A. Duran]{Angel Duran}
\address{ Applied Mathematics Department,  University of
Valladolid, 47011 Valladolid, Spain}
\email{angel@mac.uva.es}




\begin{abstract}
In this paper a semidiscrete Fourier pseudospectral method for approximating Benjamin-type equations is introduced and analyzed. A study of convergence is presented.
\end{abstract}

\maketitle

\section{Introduction}
\label{intro}
This paper is concerned with numerical approximations to Benjamin-type equations of the form
\begin{eqnarray}
u_{t}-\mathcal{L}u_{x}+f(u)_{x}=0.\label{intro1}
\end{eqnarray}
In (\ref{intro1}) $u=u(x,t)$ is a real-valued function, $\mathcal{L}$ is the linear, nonlocal, pseudodifferential operator with Fourier symbol
\begin{eqnarray}
\widehat{\mathcal{L}u}(\xi)=l(\xi)\widehat{u}(\xi)=(\delta |\xi|^{2m}-\gamma
|\xi|^{2r})\widehat{u}(\xi),\quad \xi\in\mathbb{R}, \label{intro2}
\end{eqnarray}
where $m\geq 1$ is an integer, $0\leq r<m, \gamma\geq 0, \delta>0$ and $\widehat{u}(\xi)$ denotes the Fourier transform of $u$ at $\xi$. Finally the nonlinear term $f$ is of the form
\begin{eqnarray}
f(u)=\frac{u^{q+1}}{q+1},\quad q\geq 1.\label{intro3}
\end{eqnarray}
The general form (\ref{intro1})-(\ref{intro3}) includes the case of the Benjamin equation, \cite{Benjamin1972,Benjamin1996}, and generalized versions ($m=1, r=1/2, q\geq 1$). Summarized here is some literature on (\ref{intro1}). This is mainly focused on \cite{LinaresS2005} where Linares and Scialom, based on the theory developed in \cite{KenigPV1993,KenigPV1994} for the generalized KdV and Benjamin-Ono (BO) equations, respectively, establish local and global well-posedness results for the corresponding initial-value problem for (\ref{intro1})-(\ref{intro3}). More specifically, the generalized Benjamin equation is studied first and for $u_{0}\in H^{s}(\mathbb{R}), s\geq 1$, the existence of $T=T(||u_{0}||_{H^{s}})$ and a unique solution $u\in C([0,T],H^{s}(\mathbb{R}))$ with $u(0)=u_{0}, ||u||_{L_{x}^{4},L_{T}^{\infty}}<\infty$ and
$$
\left(\int_{0}^{T}|\partial_{x} u(x,t)|^{2}dt\right)^{1/2}<\infty,\quad x\in\mathbb{R},
$$ 
are proved. Furthermore, the initial-value problem is globally well-posed in $H^{1}(\mathbb{R})$ for $q=2,3$ with no restriction on the initial data and for $q\geq 4$ when the initial data is small enough. In the general case and for $u_{0}\in H^{m}(\mathbb{R}), m>r>0, m>1$, the initial-value problem is globally well-posed in the energy space where the energy $E(u)$ given by
\begin{eqnarray}
E(u)=\int_{-\infty}^{\infty} \left(u\mathcal{L}u dx-2 F(u)\right)dx,\label{gben7}
\end{eqnarray}
is finite. Here $F^{\prime}(u)=f(u), F(0)=0$. Existence of global solutions is obtained without restriction on the initial data when $q<4m$, while for $q\geq 4m$ the initial condition must be small enough. In addition to (\ref{gben7}) and for decaying and smooth enough solutions two other quantities 
\begin{eqnarray}
I(u)=\int_{-\infty}^{\infty} u^{2} dx;\quad C(u)=\int_{-\infty}^{\infty} u dx.\label{conservedq}
\end{eqnarray}
are preserved. 

Another group of results concerns solitary-wave solutions of (\ref{intro1})-(\ref{intro3}), that is solutions of the form
$$ u=\varphi(x-c_{s}t),\quad c_{s}>0,$$ where $\varphi$ and its derivatives go to zero as $X=x-c_{s}t\rightarrow \pm \infty$. Specifically, in \cite{BonaC1998} existence and asymptotic properties of solitary waves are proved for a range of values of $\gamma$ depending on $r, m$ and the speed $c_{s}$. Stability (in the orbital sense) of the waves is studied in \cite{Angulo2009} (see also the references therein).

Typically the dynamics of these solitary waves is the main motivation to study numerically nonlinear dispersive wave equations like (\ref{intro1}), by constructing schemes approximating to corresponding initial- and periodic boundary-value problems on sufficiently large intervals. With this aim the present paper is focused on the numerical analysis. Consider, for simplicity, the initial- and periodic boundary-value problem (ipbvp) for (\ref{intro1})-(\ref{intro3}) on $[-\pi,\pi]$
\begin{eqnarray}
&&u_{t}-\mathcal{L}u_{x}+f(u)_{x}=0,\quad x\in [-\pi,\pi], t>0,\label{a1}\\
&&u(x,0)=u_{0}(x),\quad x\in [-\pi,\pi],\nonumber
\end{eqnarray}
where $u_{0}$ is $2\pi-$periodic. In this case, to our knowledge, the literature on the problem (\ref{a1}) can be summarized as follows. In \cite{Linares1999}, using the methods by Bourgain, \cite{Bourgain1993}, and Kenig et al., \cite{KenigPV1993}, for the study of the well-posedness of the ipbvp for the KdV equation, global well-posedness for the Benjamin equation is established for data in $L^{2}(\mathbb{T})$. Recently, Shi and Li, \cite{ShiL2015}, proved local well-posedness for small initial data in $H^{\mu}(\mathbb{T}), \mu\geq -1/2$, by using a different technique. Finally Cascaval, \cite{Cascaval2004}, studies the local and global well-posedness of the initial value problem (ivp) and the ipbvp for the class of nonlinear dispersive equations
$$
u_{t}-Mu_{x}+G(u)_{x}=0,
$$ where $\widehat{Mu}(\xi)=|\xi|^{2\beta}\widehat{u}(\xi), \beta\geq 1/2$ and $G$ is sufficiently smooth and satisfies
$$
\lim_{|r|\rightarrow\infty}\frac{G^{\prime}(r)}{|r|^{p}}<\infty,
$$ for some $p<4\beta$. (This includes the cases $G(u)_{x}=u^{p}u_{x}, p\geq 1$.) The author of \cite{Cascaval2004} obtains global well-posedness in $H^{\mu}(\mathbb{T})$ with $\mu=\max\{2\beta,3/2+\epsilon\}$ for some $\epsilon>0$ and in the cases:
\begin{itemize}
\item $\beta=1, \mu=2$.
\item $\beta=1/2, G^{\prime}(u)=u, \mu\in (3/2,2]$.
\item $\beta>3/2, \mu=2\beta$.
\end{itemize}
Furthermore, for $1/2<\beta\leq 3/2$ the problem is locally well-posed in $H^{\mu}(\mathbb{T})$ with $\mu=2\beta$ for $\beta>3/4$ and $\mu\in (3/2,2\beta+1]$ for $1/2<\beta\leq 3/4$. In particular, global well-posedness of the periodic gKdV equation (for $p<4$) in $H^{2}(\mathbb{T})$ is obtained. 

On the other hand, in \cite{BonaCh2013} well-posedness of the ivp of some models of the form (\ref{a1}) in relatively smooth, periodic function spaces is assumed. These spaces must have at least finite energy in the sense that the solutions
\begin{eqnarray*}
u(x,t)=\sum_{n=-\infty}^{\infty} u_{n}(t)e^{i{n x}},
\end{eqnarray*}
decomposed into its Fourier series has the property
\begin{eqnarray*}
\sum_{n=-\infty}^{\infty} \left|l\left({n}\right)\right||u_{n}(t)|^{2}<\infty,\quad t>0,
\end{eqnarray*}
where $l(\xi)$ is defined in (\ref{intro2}) as the Fourier symbol of $\mathcal{L}$. Also, for smooth solutions, periodic conditions ensure the preservation of the corresponding versions of the functionals (\ref{gben7}), (\ref{conservedq}) with the integrals on the interval $(-\pi,\pi)$, that is
\begin{eqnarray}
C_{\pi}(u)=\int_{-\pi}^{\pi} u dx;\quad 
I_{\pi}(u)=\int_{-\pi}^{\pi} u^{2} dx,\quad 
E_{\pi}(u)=\int_{-\pi}^{\pi} \left(u\mathcal{L}u dx-2 F(u)\right)dx,\label{conservedq2}
\end{eqnarray}
where $F^{\prime}(u)=f(u), F(0)=0$. 

To our knowledge, the numerical approximation to (\ref{a1}) has only been considered for the particular cases of the KdV equation and generalized KdV equation ($m=1, \gamma=0$) and for the Benjamin equation and generalized Benjamin equation ($m=1, r=1/2$). Focused on this last case, we first mention the method used in \cite{AlbertBR1999,BonaK2004} based on pseudospectral collocation in space and a second-order time-stepping code. On the other hand, the method considered in \cite{DougalisDM2012}  has the same type of spatial discretization with a third-order singly diagonally implicit Runge-Kutta scheme, combined with a projection technique to preserve invariant quantities, \cite{HairerLW2004}, as time integrator. In \cite{DougalisDM2015}, a hybrid spectral-finite element scheme along with a 2-stage Gauss-Legendre implicit Runge-Kutta method is constructed while structure preserving integrators are proposed in \cite{KinugasaMM2015}. As mentioned above, the application of all of them was, to a greater or lesser extent, related to solitary wave dynamics.

The purpose of the present paper is the introduction and analysis of a semi-discrete numerical method to approximate (\ref{a1}), based on a Fourier pseudospectral discretization. As observed before, this is, to our knowledge, the first proposal in the literature to approximate the periodic-initial value problem (\ref{a1}), (\ref{intro2}), (\ref{intro3}) in its full generality. The possible nonlocal character of the linear operator $\mathcal{L}$ in (\ref{intro2}) justifies the use of a Fourier-Galerkin discretization in space. 
{The semidiscrete scheme is studied in Section \ref{sec:1}. We first establish the existence of a unique, local in time solution of the semidiscrete problem. The existence of a global in time solution is then derived from the preservation of the quantities (\ref{conservedq2}) by the semidiscretization. The convergence of the pseudospectral discretization is then established in Theorem \ref{Theo1}. The proof is based on the introduction of a linear, intermediate problem, \cite{BakerDK1983,PelloniD2001}, whose solution approximates both those of the continuous and semidiscrete problems, with estimates that depend on the regularity of the solution of the first one}.

For real $\mu\geq 0$ and $1\leq p\leq\infty$ we denote by $W_{p}^{\mu}=W_{p}^{\mu}(-\pi,\pi)$ the real Sobolev space on $(-\pi,\pi)$ with norm $||\cdot||_{\mu,p}$. Let $H^{\mu}:=W_{2}^{\mu}$ and for $g\in H^{\mu}$ put
\begin{eqnarray*}
&&||g||_{{\mu}}=||g||_{\mu,2}=\left(\sum_{k\in\mathbb{Z}}(1+k^{2})^{\mu}|\widehat{g}(k)|^{2}\right)^{1/2},\\
&&\widehat{g}(k)=\frac{1}{2\pi}\int_{-\pi}^{\pi}e^{-ikx}g(x)dx
\end{eqnarray*}
$|\cdot|_{\infty}$ will stand for the norm in $L^{\infty}(-\pi,\pi)$. Finally, the inner product in $H^{0}=L^{2}(-\pi,\pi)$ is denoted by $(\cdot,\cdot)$, defined by
\begin{eqnarray*}
(u,v)=\int_{-\pi}^{\pi}u(x)\overline{v(x)}dx,
\end{eqnarray*}
with $||\cdot||$ standing for the corresponding norm. According to previous comments and the form of $\mathcal{L}$ in  (\ref{intro2}), we will assume that problem (\ref{a1}), (\ref{intro2}), (\ref{intro3}) is well-posed on $H^{\mu}, \mu\geq m$.
\section{Analysis of convergence of the semidiscrete scheme}
\label{sec:1}
\subsection{A Fourier pseudospectral approximation of the periodic-initial value problem}
For $N\geq 1$ integer we consider
\begin{eqnarray*}
S_{N}=span\{e^{ikx}, k\in\mathbb{Z}, -N\leq k\leq N\}.
\end{eqnarray*}
The semidiscrete Fourier-Galerkin approximation to (\ref{a1}), (\ref{intro2}), (\ref{intro3}) is defined as a real-valued map $u^{N}:[0,\infty)\rightarrow S_{N}$ such that, for all $\chi\in S_{N}$,
\begin{eqnarray}
&&(u_{t}^{N},\chi)+((-\mathcal{L} u_{x}^{N}+f(u^{N})_{x}),\chi)=0,\quad t>0,\label{a3}\\
&&u^{N}(x,0)=P_{N}u_{0}(x),\nonumber
\end{eqnarray}
where $P_{N}$ is the $L^{2}-$projection of $L^{2}$ onto $S_{N}$. For $v\in L^{2}$ we have
$$P_{N}v=\sum_{|k|\leq N}\widehat{v}_{k}e^{ikx},$$ with $\widehat{v}_{k}$ the $k-$th Fourier coefficient of $v$. Some properties of $P_{N}$ will be used throughout the paper. First, it is well known that $P_{N}$ commutes with the differential operator $\partial_{x}$. Moreover, cf. \cite{CHQZ}, given integers $0\leq j\leq\mu$ there exists a constant $C$ independent of $N$ such that for any $v\in H^{\mu}$,
\begin{eqnarray}
||v-P_{N}v||_{j}&\leq &CN^{j-\mu}||v||_{\mu},\quad \mu\geq 0,\nonumber\\ 
||v-P_{N}v||_{\infty}&\leq &CN^{1/2-\mu}||v||_{\mu},\quad \mu\geq 1.\label{21b}
\end{eqnarray}
When $\chi=e^{ikx}, |k|=0,\ldots,N$ then (\ref{a3}) becomes the initial value problem for the Fourier coefficients of $u^{N}$,
\begin{eqnarray}
\widehat{u^{N}_{t}}(k,t)&=&(ik)((\delta |k|^{2m}-\gamma |k|^{2r})\widehat{u^{N}}(k,t)-\widehat{f(u^{N})}(k,t)),\quad t>0,\label{a4}\\
&&\widehat{u^{N}}(k,0)=\widehat{u_{0}}(k),\nonumber
\end{eqnarray}
\begin{remark}
\label{remark1}
In the sequel the following properties of $f$ and $f^{\prime}$ will be used:
\begin{itemize}
\item[(i)] $f(u)-f(v)=(u-v)g(u,v,q)$ with 
$$g(u,v,q)=\frac{1}{q+1}\sum_{j=0}^{q}u^{j}v^{q-j},$$ and therefore
$$||f(u)-f(v)||\leq ||u-v|||g|_{\infty},\quad |g|_{\infty}\leq (\max\{|u|_{\infty},|v|_{\infty}\})^{q}.$$
\item[(ii)] $f^{\prime}(u)-f^{\prime}(v)=(u-v)h(u,v,q)$ with 
$$h(u,v,q)=\sum_{j=0}^{q-1}u^{j}v^{q-1-j},$$ and therefore
$$||f^{\prime}(u)-f^{\prime}(v)||\leq ||u-v|||h|_{\infty}, \quad |h|_{\infty}\leq q (\max\{|u|_{\infty},|v|_{\infty}\})^{q-1}.$$
\end{itemize}
\end{remark}
\subsection{Existence and uniqueness of solutions of the semidiscrete problem}
In particular, (i) implies that
the right hand side of (\ref{a4}) is at least locally Lipschitz continuous with respect to the $L^{2}$ norm in $S_{N}$. Then, using standard theory of ordinary differential equations, we obtain the existence of a unique, local in time solution of (\ref{a4}). Also, standard arguments prove the existence of a global in time solution if the semidiscretization preserves the $L^{2}$ norm. In our case, we have the following result.
\begin{lemma}
\label{lem21}
The solution $u^{N}$ of (\ref{a3}) satisfies, for $t>0$,
\begin{eqnarray*}
\frac{d}{dt}C_{\pi}(u^{N})=\frac{d}{dt}I_{\pi}(u^{N})=\frac{d}{dt}E_{\pi}(u^{N})=0.
\end{eqnarray*}
where $C_{\pi}, I_{\pi}$ and $E_{\pi}$ are given by (\ref{conservedq2}).
\end{lemma}

{\em Proof.} The preservation of $C_{\pi}$ is obtained directly by taking $\chi=1$ in (\ref{a3}) while if we take $\chi=u^{N}$ we have
$$(f(u^{N})_{x},u^{N})=F(u^{N}(\pi,t))-F(u^{N}(-\pi,t))=0,$$ (where $F$ is a primitive of $f$). Moreover, if $v$ is a real-valued element of $S_{N}$ it follows, since $\widehat{v}(k)=\overline{\widehat{v}(-k)}$, that
\begin{eqnarray}
(\mathcal{L} v_{x},v)=0.\label{23}
\end{eqnarray}
Therefore $(\mathcal{L} u_{x}^{N},u^{N})=0,$ which implies the preservation of $I_{\pi}$. Finally, with $\chi=P_{N}f(u^{N})-\mathcal{L}u^{N}$ we have
\begin{eqnarray*}
(u_{t}^{N},\chi)&=&\int_{-\pi}^{\pi} u_{t}^{N}(f(u^{N})-\mathcal{L}u^{N})dx\\
&=&-\frac{1}{2}
\frac{d}{dt}\int_{-\pi}^{\pi} \left(u^{N}\mathcal{L}u^{N}-\frac{2}{(q+1)(q+2)}(u^{N})^{q+2}\right)dx\\
&=&-\frac{1}{2}\frac{d}{dt}E_{\pi}(u^{N}(\cdot,t)),
\end{eqnarray*}
and because of periodicity and properties of $P_{N}$, we have
\begin{eqnarray*}
(-\mathcal{L} u_{x}^{N}+f(u^{N})_{x},\chi)&=&(-\mathcal{L} u_{x}^{N}+f(u^{N})_{x},P_{N}f(u^{N})-\mathcal{L}u^{N})\\
&=&(f(u^{N})_{x},f(u^{N}))-(-\mathcal{L} u_{x}^{N},f(u^{N}))\\
&&
-(f(u^{N})_{x},\mathcal{L}u^{N})+(\mathcal{L} u_{x}^{N},\mathcal{L}u^{N})
=0.\quad 
\qed
\end{eqnarray*}
\subsection{Analysis of convergence}
In order to study the convergence of the semidiscrete scheme, we consider the intermediate problem of searching for $w^{N}\in S_{N}$ such that for all $\chi\in S_{N}$ and any $t^{*}>0$, \cite{BakerDK1983,PelloniD2001},
\begin{eqnarray}
&&(w_{t}^{N},\chi)+(-\mathcal{L} w_{x}^{N}+f^{\prime}(u)w^{N}_{x},\chi)=0,0\leq t\leq t^{*},\label{a5}\\
&&w^{N}(0)=P_{N}u_{0},\nonumber
\end{eqnarray}
where $u$ is a solution of (\ref{a1}). 

\begin{lemma}
\label{lem22}
Let $u$, the solution of (\ref{a1}), belong to $H^{\mu}$ for $\mu>3/2$. Then the ivp (\ref{a5}) has a unique solution $w^{N}$ which satisfies
\begin{eqnarray}
\max_{0\leq t\leq t^{*}}|w^{N}|_{\infty}\leq C,\label{25}
\end{eqnarray}
and
\begin{eqnarray}
\max_{0\leq t\leq t^{*}}||u-w^{N}||\leq C N^{1-\mu},\label{26}
\end{eqnarray}
for some constant $C=C(u,t^{*})$. In addition, if $\mu\geq 5/2$ then
\begin{eqnarray}
\max_{0\leq t\leq t^{*}}||w^{N}||_{1,\infty}\leq C,\label{27}
\end{eqnarray}
where again $C$ is a constant depending on $u$ and $t^{*}$ only.
\end{lemma}

{\em Proof.} As before, standard theory ensures existence and uniqueness of a local solution for (\ref{a5}). Furthermore, while $w^{N}$ exists, putting $\chi=w^{N}$ in (\ref{a5}) gives, in view of (\ref{23}), that
$$
\frac{d}{dt}\frac{||w^{N}||^{2}}{2}+(f^{\prime}(u)w^{N}_{x},w^{N})=0,
$$ that is
\begin{eqnarray*}
\frac{d}{dt}\frac{||w^{N}||^{2}}{2}=\frac{1}{2}(\partial_{x}(f^{\prime}(u)),(w^{N})^{2})
\leq \frac{1}{2}|\partial_{x}(f^{\prime}(u))|_{\infty}||w^{N}||^{2}.
\end{eqnarray*}
Now
$$
|\partial_{x}(f^{\prime}(u))|_{\infty}=|u^{q-1}u_{x}|_{\infty},
$$ which is bounded for all $t>0$ if $\mu> 3/2$. Then by Gronwall's lemma, there is a constant $C$ such that
$$||w^{N}||\leq Ce^{Ct},$$
and we can extend the local solution to a solution on $[0,t^{*}]$.

Now we define $\rho^{N}=P_{N}u-w^{N}$ and estimate the difference
\begin{eqnarray*}
u-w^{N}=u-P_{N}u+\rho^{N}.
\end{eqnarray*}
In view of (\ref{a5}) $\rho^{N}$ satisfies, for all $\chi\in S_{N}$,
\begin{eqnarray}
&&(\rho_{t}^{N},\chi)+(-\mathcal{L} \rho_{x}^{N}+f(u)_{x}-f^{\prime}(u)w^{N}_{x}),\chi)=0,\quad 0\leq t\leq t^{*},\label{a8}\\
&&\rho^{N}(0)=0,\nonumber
\end{eqnarray}
Note first that for $\chi\in S_{N}$
\begin{eqnarray*}
(f(u)_{x}-f^{\prime}(u)w^{N}_{x}),\chi)
=(f^{\prime}(u)[(u-P_{N}(u))_{x}+\rho^{N}_{x}],\chi).
\end{eqnarray*}
Therefore, if we take $\chi=\rho^{N}$ in (\ref{a8}) we have
\begin{eqnarray*}
\frac{1}{2}\frac{d}{dt}||\rho^{N}||^{2}+(f^{\prime}(u)(u-P_{N}(u))_{x},\rho^{N})+(f^{\prime}(u)\rho^{N}_{x},\rho^{N})=0,
\end{eqnarray*}
and by periodicity
$$
(f^{\prime}(u)\rho^{N}_{x},\rho^{N})=-\frac{1}{2}(\partial_{x}f^{\prime}(u),(\rho^{N})^{2}).
$$ Hence
\begin{eqnarray*}
\frac{1}{2}\frac{d}{dt}||\rho^{N}||^{2}&=&-(f^{\prime}(u)(u-P_{N}(u))_{x},\rho^{N})+
\frac{1}{2}(\partial_{x}f^{\prime}(u),(\rho^{N})^{2})\\
&\leq&|f^{\prime}(u)|_{\infty}||(u-P_{N}(u))_{x}||||\rho^{N}||+\frac{1}{2}|\partial_{x}f^{\prime}(u)|_{\infty}||\rho^{N}||^{2}\\
&\leq & \frac{C}{N^{\mu-1}}||\rho^{N}||+C||\rho^{N}||^{2}\\
&\leq & \left(\frac{C}{N^{\mu-1}}\right)^{2}+||\rho^{N}||^{2}.
\end{eqnarray*}
Finally, the initial condition and Gronwall's lemma lead to
\begin{eqnarray}
\max_{0\leq t\leq t^{*}}||\rho^{N}||\leq \frac{C}{N^{\mu-1}},\label{a9}
\end{eqnarray}
and, consequently, to (\ref{26}).

Note now that if $0\leq t\leq t^{*}$
\begin{eqnarray*}
|w^{N}(t)|_{\infty}&\leq &|u(t)|_{\infty}+ |w^{N}(t)-P_{N}u(t)|_{\infty}+
|P_{N}u(t)-u(t)|_{\infty}\\
&=&|u(t)|_{\infty}+|\rho^{N}(t)|_{\infty}+
|P_{N}u(t)-u(t)|_{\infty}
\end{eqnarray*}
Recall now that the following inverse inequalities hold in $S_{N}$. Given $0\leq s\leq r$, there exists a constant $C_{0}$ such that
\begin{eqnarray}
||\psi||_{r}\leq C_{0}N^{r-s}||\psi||_{s},\quad 
||\psi||_{r,\infty}\leq C_{0}N^{1/2+r-s}||\psi||_{s},\label{inv_ineq}
\end{eqnarray}
for all $\psi\in S_{N}$. Then (\ref{21b}) and (\ref{a9}) give
\begin{eqnarray*}
|w^{N}(t)|_{\infty}\leq C+\frac{C}{N^{\mu-3/2}}+\frac{C}{N^{\mu-1/2}}\leq \tilde{C},\label{a14}
\end{eqnarray*}
for some constant $C$ and (\ref{25}) follows. If now $\mu\geq 5/2$ we have
\begin{eqnarray*}
|w_{x}^{N}(t)|_{\infty}&\leq&|u_{x}(t)|_{\infty}+|\rho_{x}^{N}(t)|_{\infty}+
|\partial_{x}(P_{N}u(t)-u(t))|_{\infty}\nonumber\\
&
\leq &C+\frac{C}{N^{\mu-5/2}}+\frac{C}{N^{\mu-3/2}}\leq \tilde{C},\label{a15}
\end{eqnarray*}
Hence (\ref{27}) follows. $\qed$

We now proceed to study the convergence of the semidiscrete scheme (\ref{a3}). 
\begin{theorem}
\label{Theo1}
Let $u^{N}$ be the solution of (\ref{a3}) and suppose that $u$, the solution of (\ref{a1}), belongs to $H^{\mu}, \mu\geq 5/2$. Then given $0\leq t^{*}<\infty$ there exists a constant $C=C(u,t^{*})$ such that
\begin{eqnarray}
\max_{0\leq t\leq t^{*}}||u-u^{N}||\leq \frac{C}{N^{\mu-1}}.\label{210}
\end{eqnarray}
\end{theorem}

{\em Proof.} We have already established the existence of $u^{N}$, the solution of the initial-value problem (\ref{a3}) on any temporal interval $[0,t^{*}]$. Let $w^{N}$ be the solution of the intermediate problem (\ref{a5}) and $e^{N}=w^{N}-u^{N}$. Since $u-u^{N}=u-w^{N}+e^{N}$, in view of (\ref{26}) we only need to estimate $e^{N}$. This satisfies, for all $\chi\in S_{N}$
\begin{eqnarray}
&&(e_{t}^{N},\chi)+(-\mathcal{L} e_{x}^{N}-(f(u^{N})_{x}-f^{\prime}(u)w^{N}_{x}),\chi)=0,\quad 0\leq t\leq t^{*},\label{a10}\\
&&e^{N}(0)=0,\nonumber
\end{eqnarray}
Since $f(u^{N})_{x}-f^{\prime}(u)w^{N}_{x}=(f(u^{N})_{x}-f(w^{N})_{x})+(f(w^{N})_{x}-f^{\prime}(u)w^{N}_{x})$, taking $\chi=e^{N}$ in (\ref{a10}) gives, in view of (\ref{23}), that for $0\leq t\leq t^{*}$
\begin{eqnarray}
(e_{t}^{N},e^{N})=(f(u^{N})_{x}-f(w^{N})_{x},e^{N})+((f^{\prime}(w^{N})-f^{\prime}(u))w^{N}_{x}),e^{N}).\label{a12}
\end{eqnarray}
In order to estimate the second term in the right hand side of (\ref{a12}), using the inclusion of Remark \ref{remark1} (ii) we observe that
\begin{eqnarray*}
||f^{\prime}(w^{N})-f^{\prime}(u)||\leq q
||w^{N}-u||\left(\max(|w^{N}|_{\infty},|u|_{\infty})\right)^{q-1}.
\end{eqnarray*}
Therefore, by (\ref{26}) , (\ref{27})
\begin{eqnarray}
|((f^{\prime}(w^{N})-f^{\prime}(u))w^{N}_{x},e^{N})|&\leq &C||f^{\prime}(w^{N})-f^{\prime}(u)|| ||w_{x}^{N}||_{\infty}||e^{N}||\nonumber\\
&\leq &\frac{C}{N^{\mu-1}} ||e^{N}||,\label{a16b}
\end{eqnarray}
for some constant $C=C(u,t^{*})$.

In order to estimate the first term in the right hand side of (\ref{a12}), we make use of formulas (3.10)-(3.13) of \cite{BonaDKM1995} to write
\begin{eqnarray*}
f(u^{N})_{x}-f(w^{N})_{x}=f(w^{N}-e^{N})_{x}-f(w^{N})_{x}=
f(-e^{N})_{x}+R(w^{N},-e^{N}),
\end{eqnarray*}
where
\begin{eqnarray*}
R(v,w)=\frac{1}{q+1}\partial_{x}\left(\sum_{j=1}^{q}\begin{pmatrix}q+1\\j\end{pmatrix}v^{q+1-j}w^{j}\right).\label{aR}
\end{eqnarray*}
Then, by periodicity
\begin{eqnarray}
(f(u^{N})_{x}-f(w^{N})_{x},e^{N})=((R(w^{N},-e^{N}),e^{N}).\label{a17}
\end{eqnarray}
Hence, using the estimate (3.13) of \cite{BonaDKM1995} we have
\begin{eqnarray}
|(R(w^{N},-e^{N}),e^{N})|\leq C_{q}\max_{1\leq m\leq q}||w^{N}||_{1,\infty}^{m}
\sum_{j=1}^{q}\int_{-\pi}^{\pi}|e^{N}|^{j+1}dx,\label{a18}
\end{eqnarray}
and we observe that
\begin{eqnarray*}
\sum_{j=1}^{q}\int_{-\pi}^{\pi}|e^{N}|^{j+1}dx\leq \max_{1\leq j\leq q}|e^{N}|_{\infty}^{j-1}||e^{N}||^{2}.
\end{eqnarray*}
Since $e^{N}(0)=0$, there exists by continuity a maximal temporal value $t_{N}>0$ (assume without loss of generality that $t_{N}<t^{*}$) such that
\begin{eqnarray}
|e^{N}|_{\infty}\leq 1,\quad 0\leq t\leq t_{N}.\label{217}
\end{eqnarray}
For $t\leq t_{N}$ therefore, (\ref{a17}), (\ref{a18}), (\ref{27}) and (\ref{217}) yield
\begin{eqnarray}
|(f(u^{N})_{x}-f(w^{N})_{x},e^{N})|\leq C||e^{N}||^{2},\label{218}
\end{eqnarray}
for some constant $C=C(u,t^{*})$. Hence by (\ref{a12}), (\ref{a16b}), and (\ref{218}) we have for $t\in [0,t_{N}]$ 
\begin{eqnarray*}
\frac{1}{2}\frac{d}{dt}||e^{N}||^{2}\leq C||e^{N}||^{2}+\frac{C}{N^{\mu-1}} ||e^{N}||
\leq  C||e^{N}||^{2}+\left(\frac{C}{N^{\mu-1}}\right)^{2}.
\end{eqnarray*}
Since $e^{N}(0)=0$, Gronwall's lemma implies
\begin{eqnarray}
\max_{0\leq t\leq t_{N}}||e^{N}||\leq \frac{C}{N^{\mu-1}},\label{a20}
\end{eqnarray}
for some constant $C=C(u,t^{*})$. Hence, by the inverse properties (\ref{inv_ineq}) of $S_{N}$ we obtain that 
$$\max_{0\leq t\leq t_{N}}|e^{N}|_{\infty}\leq CN^{3/2-\mu},$$ which contradicts the maximality of $t_{N}$ in (\ref{217}) if $N$ is taken sufficiently large. We conclude that $t_{N}=t^{*}$ and (\ref{210}) follows from (\ref{a20}). $\qed$

\begin{remark}
Under the hypothesis of Theorem \ref{Theo1} we have for $0\leq t\leq t^{*}$, using (\ref{27}), inverse properties of $S_{N}$, (\ref{26}) and (\ref{210})
\begin{eqnarray*}
||u^{N}||_{1,\infty}&\leq &||u^{N}-w^{N}||_{1,\infty}+||w^{N}||_{1,\infty}\\
&\leq & CN^{3/2}||u^{N}-w^{N}||+C\\
&\leq & CN^{3/2}\left(||u^{N}-u||+||u-w^{N}||\right)+C\leq CN^{5/2-\mu}+C,
\end{eqnarray*}
i.e.
\begin{eqnarray*}
||u^{N}||_{1,\infty}\leq B,\label{220}
\end{eqnarray*}
for some constant $B=B(u,t^{*})$. 
\end{remark}


\section*{acknowledgements}
This work was supported by Spanish MINECO under Research Grants MTM-54710-P and TEC2015-69665-R and by JCYL under Research Grant VA041P17



\end{document}